\documentclass[11pt]{article}
\usepackage{latexsym,amssymb,amsmath}
\def\lra{\longrightarrow}
\def\ra{\rightarrow}
\def\CC{\mathbb C}
\def\PP{\mathbb P}
\def\QQ{\mathbb Q}

\def\cC{\mathcal C}
\def\cF{\mathcal F}
\def\cH{\mathcal H}
\def\cI{\mathcal I}
\def\cJ{\mathcal J}
\def\cL{\mathcal L}
\def\cO{\mathcal O}
\def\cM{\mathcal M}
\def\cX{\mathcal X}

\def\fg{\mathfrak g}
\def\fsl{\mathfrak sl}
\def\s{\sigma}
\def\a{\alpha}
\def\op{\oplus}

\newcommand\proof{{\noindent {\em Proof}.}\hspace{2mm}}
\newcommand\qed{{\hfill\hfill $\Box$}}
\newtheorem{theo}{Theorem}
\newtheorem{coro}[theo]{Corollary}
\newtheorem{lemm}[theo]{Lemma}
\newtheorem{prop}[theo]{Proposition}

\begin{document}
\title{Prime Fano threefolds and integrable systems}
\author{A. Iliev, L. Manivel}
\date{}
\maketitle
\footnotetext[1]{\ MSC 14J28, 14J45, 14K30, 37J99}
\footnotetext[2]{\ The 1-st author is partially supported by grant MI1503/2005
                  of the Bulgarian Foundation for Scientific Research}
\begin{abstract}
For a general K3 surface $S$ of genus $g$, $2 \le g \le 10$,
we prove that the intermediate Jacobians of the family of prime
Fano threefolds of genus $g$ containing $S$ as a hyperplane section,
form generically an algebraic completely integrable system.
\end{abstract}

\section{Introduction}

The punctual Hilbert schemes $Hilb_dS$ of length $d$ subschemes
of a smooth K3 surface $S$ are examples of smooth complex 
symplectic varieties \cite{beau1}. The K3 surfaces are the only 
smooth compact complex symplectic varieties of dimension two.
Any curve $C$ on a K3 surface is automatically a Lagrangian subvariety
of $S$, and if $C \subset S$ is a smooth curve of genus $g$
then the relative Jacobian $\cJ$
is a Lagrangian fibration (equivalently, an algebraically completely
integrable system, see \cite{DM}, 2.3)
over an open subset of the complete linear system
$|C| \simeq \PP^g$, see \cite{beau2}.
In fact, if all the curves in $|C|$ are reduced,   then
the compactified relative Jacobian $\bar\cJ \rightarrow |C| \simeq \PP^g$
is a smooth symplectic variety and still remains a Lagrangian fibration,
see also \S 4 in \cite{Bot}.
Recently Sawon, by using the techniques of the Fourier-Mukai functors,
managed to show that for any $k\ge 2$, and for
the general K3 surface $S$ with primitive polarization of degree
$d=k^2(2g-2)$, the Hilbert scheme $Hilb_gS$ is locally isomorphic
to a compactified relative Jacobian $\bar\cJ \rightarrow \PP^g$ as
above, thus proving the existence of a Lagrangian fibration on
$Hilb_gS$ over $\PP^g$, see \cite{Saw}. The fibrations
of Beauville and Sawon are fibered generically by Jacobians of
curves, and one may ask how to find other types of Lagrangian
fibrations.

A way towards new examples is the following criterion of integrability,
due to Donagi and Markman and generalizing the above construction
of Beauville (see 8.1-8.2 in \cite{DM}):

\medskip\noindent
(DM)
{\it Let $M$ be a smooth complex symplectic variety, and let
$Z \subset M$ be a smooth Lagrangian subvariety.
Then the relative Picard bundle $Pic^o$ is a Lagrangian fibration
(an algebraically completely integrable system) over an open set of the
base of deformations of $Z$ in $M$.}

\medskip
In general, the criterion (DM) raises the problem of finding examples 
of complex symplectic varieties $M$ (different from K3 surfaces) 
and Lagrangian subvarieties $Z$ with well understood Picard variety 
and deformations in $M$. 

The first  higher dimensional symplectic variety 
that has been exploited in this 
context is the Beauville-Donagi variety -- the Fano fourfold 
$F(W)$ of lines on the general smooth cubic fourfold $W$, see \cite{BD}.  
This is performed in two different setups, corresponding to 
two basic families of Lagrangian surfaces in $F(W)$ arising 
from the geometry of cubic hypersurfaces of dimensions $3$ and $5$.
%
In the first case, as observed by Voisin, 
the Fano families of lines on the general cubic threefolds contained 
in $W$ as hyperplane sections are smooth Lagrangian surfaces in $F(W)$, 
see Exemple 3 in \cite{voisin}. Applying (DM), one obtains 
(see Example 8.22 of \cite{DM}):
{\it The relative intermediate Jacobian over the family of 
smooth cubic threefolds, contained as hyperplane sections of 
$W$, is an algebraic integrable system}. 
In the second case, by our recent observation (see \cite{IMn2}) 
the Fano families of planes on the general cubic fivefolds 
containing $W$ as a hyperplane section are embedded 
by the intersection map as smooth Lagrangian surfaces in $F(W)$. 
This yields the second application of (DM), which is in some 
sense   symmetric to the first: 
{\it The relative intermediate Jacobian over the family of 
smooth cubic fivefolds, containing $W$ as a hyperplane section, 
is an algebraic integrable system}, ibid. 

\medskip
In the present paper we consider, instead of the Beauville-Donagi
variety, the Hilbert powers $Hilb_dS$ of smooth K3 surfaces $S$
contained in prime Fano threefolds.

Recall that a prime Fano threefold is a smooth compact complex variety $X$
of dimension $3$ for which the anticanonical divisor $H = -K_X$ is
the ample generator of $Pic(X)$.
For such $X$, the number $g = -K_X^3/2 +1$
is always a positive integer, called the genus of $X = X_{2g-2}$.
Prime Fano threefolds exist only for
$2 \le g \le 12,\; g \not=11$ (see \cite {IP}). 
A general member of the anticanonical linear system $|-K_X|$ of
a prime Fano threefold $X = X_{2g-2}$, is a K3 surface $S = S_{2g-2}$
of genus $g$. Mukai proved that the converse is also true:
the general K3 surface of genus $g$, $2 \le g \le 12,\; g \not= 11$ 
can be embedded as an anticanonical divisor in a prime Fano threefold
of genus $g$, see \cite{mukai1} or \cite{beau3} for a more
precise statement. In this situation takes place the following 
observation,  
due originally to Thomas (see Theorem 1.6.1 in \cite{Tho}):  

\medskip
{\it Let $S = S_{2g-2}$ be a K3 surface in a prime Fano threefold
$X = X_{2g-2}$ of genus $g$.
Then for any integer $d \ge 1$, any component $\cH$ of the Hilbert
scheme of smooth curves of degree $d$ on $X$ is sent by the
intersection map $j$ with $S$ to a Lagrangian
subvariety $j(\cH)$ of $Hilb_dS$},

\medskip\noindent
see Lemma 1. Since our study will require a closer inspection of  the 
intersection map $j$,  in \S 2 we write down a detailed 
proof of Lemma 1. In particular, we state a criterion for $j$ 
to be an immersion (see Corollary \ref{criterion}).

In \S 3 and \S 4.1 we study the scheme $F(X)$
of conics on a prime Fano threefold $X = X_{2g-2}$.
Notice that together with the family of lines, this is one of the two
best studied families of curves on prime Fano threefolds;  
see Ch.4 of \cite{IP} for a survey of the basic results about lines
and conics on Fano threefolds.
In \S 3.1 we collect the most important, at least from our perspective,
 of these results,  beginning with the fact that for $X$ 
general, $F(X)$ is a smooth irreducible surface,
{\it the Fano surface} of $X$.
In \S 3.2 we study the intersection map $j: F(X) \rightarrow Hilb_2 S$ 
for the general pair 
$(X,S)$ of a prime Fano threefold $X$ and a K3 surface $S \in |-K_X|$.  
By Lemma 1, 
$j(F(X))$ is a Lagrangian surface in the symplectic fourfold
$Hilb_2S$.  

In order to apply the integrability criterion (DM), we need to know 
whether $j(F(X))$ is smooth. 
By the Proof of Proposition \ref{immers}, for $g \ge 7$ the map 
$j$ is injective; in particular the isomorphic image $j(F(X))$ 
of the smooth Fano surface $F(X)$ is a smooth Lagrangian subvariety 
of $Hilb_2S$. So for $g \ge 7$ we can directly apply the criterion 
(DM) to get an integrable system.
However for $g \le 6$ the image $j(F(X))$ is not smooth,  
even for the general $X$, and the criterion (DM) 
cannot be directly applied to pairs $(Z,M) = (j(F(X)), Hilb_2S)$.  
More concretely, for any $g \le 6$ the map $j$ is almost an 
isomorphism from $F(X)$ onto its image, except for a finite 
number of multiple points (presumably simple double points).

Nevertheless, by Proposition \ref{immers}, 
for the general pair $(X,S)$, of a prime Fano 
threefold $X$ of any possible genus $g$ and a K3 surface $S \in |-K_X|$
the intersection map $j: F(X) \rightarrow Hilb_2S$ 
is in fact an {\it immersion}.
In this situation, due to results of Ran and Voisin, 
we still can apply the criterion (DM), by replacing 
in its conditions 
the smooth Lagrangian subvariety by a Lagrangian immersion 
of a smooth variety, see \S 4.2. 
This gives rise from $j(F(X)) \subset Hilb_2S$ 
to an algebraic integrable system
(i.e. a Lagrangian fibration) -- 
the relative Picard $Pic^o \rightarrow \cH$
over the base  $\cH$ of the deformation space 
of $j$ (with a fixed target), see \S 4.2.  

It turns out that we can exactly describe the fibers
of this Lagrangian fibration.  We make this in two steps. 
In general, the fibers of such a relative Picard 
are the Picard schemes of certain smooth surfaces $F_s$, 
obtained as deformation of a fixed $F(X)=F_o$.  
In our particular case, Proposition 9
shows that all these surfaces $F_s$
are in fact Fano surfaces of Fano threefolds $X_s$
containing the K3 surface $S$,  
at least in a neighbourhood of $F(X)$.
Next, due to a less known but general fact 
about prime Fano threefolds, the Abel-Jacobi map 
from the Albanese variety $Alb\, F(X)$ of $F(X)$ 
to the intermediate Jacobian $J(X)$
of $X$ is an isomorphism, at least for the general 
prime $X$ of any genus $g$, see Theorem \ref{AJ}. 
Now by Theorem \ref{AJ}, and by the fact that the intermediate 
Jacobians $J(X_s)$ are principally polarized abelian varieties,
$J(X_s) \simeq Alb\, F(X_s) \simeq Pic^o F(X_s)$. 
This yields the main result of the paper, 
Theorem \ref{main}:

\medskip

{\it For a general K3 surface $S$ of genus $2 \le g \le 10$,
the relative intermediate Jacobian is generically an algebraic
integrable system over the family of prime Fano threefolds
containing $S$.}

\medskip
There is a strong analogy with 
the fact that the relative intermediate Jacobian is generically a
completely integrable system over a complete deformation of
gauged Calabi-Yau threefolds (\cite{DM}, Theorem 7.7).
Here we considered families of {\it log Calabi-Yau varieties with
a fixed boundary}, that is, families of pairs $(X,S)$ with $X$ a
(Fano) threefold and $S$ a divisor such that $K_X+S=0$
(the gauge disappears).
The boundary $S$ is our fixed K3 surface.
It would be interesting to extend our main result to more general
families of log Calabi-Yau threefolds.

\section{Lagrangian images of families of curves on Fano threefolds}

Recall from \cite{mukai1} and \cite{beau3} that for any $2 \le g \le
12, g \not=11$,
the general K3 surface $S = S_{2g-2}$ of genus $g$ can be represented
as an anticanonical divisor in a Fano threefold $X = X_{2g-2}$ of genus $g$.
For $g\ge 3$, the anticanonical bundle $-K_X$ is very ample,
so $S$ is just a hyperplane section
in the anticanonical embedding of $X$.

\begin{lemm}
Let $X = X_{2g-2}$ be a Fano threefold of genus $g$,
$2 \le g \le 12, g \not=11$, and let $S$ be a general member of the
anticanonical system $|-K_X|$.
Let $\cH$ be an irreducible component of the Hilbert scheme of curves on $X$,
whose general member is a smooth curve of anticanonical degree $d$,
and consider the intersection map
$$
j : \cH \ra Hilb_dS, \ \ \ j:C \mapsto C\cap S.
$$
Then the image $j(\cH) \subset Hilb_dS$ is a Lagrangian subvariety of
$Hilb_dS$.
\end{lemm}

\proof
Let the curve $C \subset X$ be a general element of $\cH$.
Since the anticanonical divisor is base-point free, and
 $S \in |-K_X|$ is also general,
the intersection $C\cap S=\{x_1,\ldots ,x_d\}$ is a reduced
zero-scheme of length $d$. Consider the differential
$$d_{\scriptscriptstyle C}j : T_C\cH=H^0(N_{C/X})\lra T_{C\cap S}Hilb_dS
= T_{x_1}S\op\cdots\op T_{x_d}S.$$

The Hilbert scheme $Hilb_dS$ has a symplectic structure
$$\wedge^2T_{C\cap S}Hilb_2S\lra \wedge^2T_{x_1}S\op\cdots
\op \wedge^2T_{x_d}S
\stackrel{\s_{x_1},\ldots ,\s_{x_d}}{\lra}\CC\op\cdots \op\CC
\stackrel{+}{\lra}\CC$$
defined by the choice of a non zero section
$\s\in H^0(S,K_S) = H^0(S,\Omega^2_S) \simeq \CC$.

Denote by $j_x$ the skew-symmetric form defined by the composition
$$
\wedge^2H^0(N_{C/X})\ra
\wedge^2T_xS\ \stackrel{\s_x}{\ra}\ \CC.$$

We shall see first that $j_{x_1}+\cdots +j_{x_d}=0$,
i.e. that $j(\cH)$ is an isotropic
subvariety of $Hilb_dS$ at the point $j(C) = x_1+\cdots+x_d$
with respect to the symplectic structure on $Hilb_dS$.
For this, consider the diagram
$$\begin{array}{cccccc}
 & & 0 & & &\\
 & & \uparrow & & &\\
 & & \cO_X(S) & & &\\
 & & \uparrow & & &\\
0\rightarrow T_xC & \rightarrow & T_xX & \rightarrow &
N_{C/X,x}&\rightarrow 0 \\
 & & \uparrow & & \uparrow & \\
 & & T_xS & \leftarrow & H^0(N_{C/X}) &\\
 & & \uparrow & & & \\
 & & 0 & &
\end{array}$$
The vertical sequence gives a natural isomorphism
$K_X(S)_{|S}=K_S$ (residue). The horizontal sequence gives
$K_{X|C}=K_C\otimes (\det N_{C/X})^*.$
So $\s$ lifts uniquely to a section of $K_X(S)$, a meromorphic
3-form on $X$ with poles on $S$. By restriction to $C$, we get a
section $\tilde{\s}$ of $K_C(C\cap S)\otimes (\det N_{C/X})^*$.
Let $\alpha,\beta$ be global sections of $N_{C/X}$. Then
$\alpha\wedge\beta$ is a global section of $\det N_{C/X}$, and
its product with $\tilde{\s}$  gives a section $(\s,\alpha\wedge\beta)$
of $K_C(C\cap S)$,
a meromorphic form on $C$ with poles at $x_1,\ldots,x_d$.
Now, for any $x \in S$ the following identity takes place

\medskip

\noindent
(*)\hspace{4.5cm} $j_x(\alpha,\beta)=res_x(\s,\alpha\wedge\beta).$

\medskip\noindent

To verify (*), choose local coordinates $u,v,w$ on $X$ around $x$,
such that $u,v$ are local equations of $C$ and $w$ a local equation
of $S$. Our generator $\s$ of $K_S$ can be written locally as $f
du\wedge dv$ for $f$ some regular function on $S$,  and it
is the residue of the form $\tilde{\s}$ on $X$ which is given
locally by $$ \tilde{\s}=Fdu\wedge dv\wedge\frac{dw}{w}, \qquad
F_{|S}=f.$$
Since $\partial/\partial u$ and $\partial/\partial v$ are
local generators of $N_{C/X}$, we write locally $\a\wedge\beta$ as
$a\partial/\partial u\wedge \partial/\partial v$. The result
of the contraction with the restriction of $\tilde{\s}$ to $C$
is $(\s,\alpha\wedge\beta)=F_{|C}a\frac{dw}{w}$, whose residue at $x$ is
$res_x(\s,\alpha\wedge\beta)=F(x)a(x)=f(x)a(x)$. On the other hand,
$j_x(\a,\beta)$
is obtained by evaluating $\s_x$ on $a(x)\partial/\partial u\wedge
\partial/\partial v$ considered as a vector of $\wedge^2T_xS$.
This gives  $j_x(\a,\beta)=f(x)a(x)$, as claimed.

\medskip

Next, by (*) and the residue theorem,
 $$j_{x_1}(\alpha,\beta)+\cdots +j_{x_d}(\alpha,\beta)
=res_{x_1}(\s,\alpha\wedge\beta)+\cdots +res_{x_d}(\s,\alpha\wedge\beta)=0,$$
which proves that $j(\cH)$ is an isotropic subvariety of
$Hilb_dS$ at $x_1+\cdots+x_d$.
\medskip

Since $x_1+\cdot+x_d$ as above is a general point of $j(\cH)$,
then $j(\cH)$ is an isotropic subvariety of $Hilb_dS$;
and in order to prove that $j(\cH)$ is a Lagrangian subvariety
of $Hilb_dS$ it remains to see that $j(\cH)$ has dimension $d$.

For this, consider the differential of $j$ at a general
member $C$ of the family $\cH$.

Let $I_C$ and $N = N_{C/X}$ denote the ideal sheaf
and the normal sheaf of $C$ in $X$.
Since the general $C \in \cH$ is smooth, in particular $C$ is a locally complete
intersection in $X$, then $N_{C/X}$ is locally free.
Let
$$T_C\cH=H^0(N)=Hom(I_C/(I_C)^2,\cO_C)\ \stackrel{d_{\scriptscriptstyle C}j}{\lra} \
T_{C\cap S}Hilb_dS=Hom(I^S_{C\cap S}/(I^S_{C\cap S})^2,\cO_{C\cap S})$$
be the natural map deduced from the identity $I^S_{C\cap S}=
(I_C+I_S)/I_S$ for the ideal sheaf of $C\cap S$ in $S$, and let
$$T_{C\cap S}Hilb_dS=Hom(I^S_{C\cap S}/(I^S_{C\cap S})^2,\cO_{C\cap S})
\stackrel{\theta}{\lra} Hom(I_C/(I_C)^2,\cO_{C\cap S})=
N\otimes_{\cO_X}\cO_{C\cap S}$$
be the map deduced from the restriction $I_C\ra I^S_{C\cap S}$.
Since this restriction map is surjective, $\theta$ is injective.
Moreover, the composition
$$H^0(N)\ \stackrel{\theta\circ d_{\scriptscriptstyle C}j}{\lra}\ N\otimes_{\cO_X}\cO_{C\cap S}$$
is just the evaluation map $ev_{\scriptscriptstyle C\cap S}$ on the punctual scheme
$C\cap S$.
In particular, the kernel of $d_{\scriptscriptstyle C}j$ coincides with the kernel of
$ev_{\scriptscriptstyle C\cap S}$, which can easily be computed from the exact sequence
$0\ra\cO_C(-S)\ra\cO_C\ra\cO_{C\cap S}\ra 0$, twisted by the locally
free sheaf $N$. We get
$$Ker\,d_{\scriptscriptstyle C}j=H^0(N(-S))=H^0(N\otimes K_{X|C})=H^0(N^*\otimes K_C)
\simeq H^1(N)^*.$$
This implies that the relative dimension of $j$ is $h^1(N)$.
Then, if $g$ denotes the arithmetic genus of $C$,
the Riemann-Roch yields
$$\dim j(\cH) = h^0(N)-h^1(N) =\chi (N)
= deg N+2(1-g)= -K_X.C =d.$$
This concludes the proof. \qed

\medskip

\noindent {\bf Example}.
Let $X$ be a prime Fano threefold of genus $g\ge 3$. Then $-K_X$
is very ample and defines an embedding $X\subset\PP^{g+1}$.
Consider the family $\cH_{can}$ of canonical
curves on $X$, that is, codimension two linear sections of $X$.
Then $\dim\cH_{can}=\dim G(2,g+2) = 2g$, while these curves have
anticanonical degree $-K_X^3=2g-2 < \dim\cH_{can}$.
The image of the map
 $$j: \cH_{can}\lra Hilb_{2g-2}S$$
is the family of canonical punctual schemes on $S$, defined again
by codimension two linear sections of $S$. Its dimension is
$\dim j(\cH_{can})=\dim G(2,g+1) = 2g-2$, as expected, and $j$
is a $\PP^2$-fibration. Note that a canonical curve $C$ has normal
bundle $N=\cO_C(1)\oplus\cO_C(1) =K_C\oplus K_C$, so clearly
$h^1(N)=2$. The fact that $j(\cH_{can})$ is Lagrangian already
appears in \cite{voisin}, Exemple 3.

\medskip
The proof of the previous statement also provides a criterion for the
map $j$ to be immersive.

\begin{coro}\label{criterion}
Let $\cH$ be as above, and let $C$ be a curve of the family parameterized
by $\cH$. Suppose that $C$ is a locally complete intersection in $X$,
 with normal bundle $N$.
Then $j$ is an immersion at $C$ if and only if the evaluation
map $$ev_{\scriptscriptstyle C\cap S}:H^0(N)\lra N\otimes_{\cO_X}\cO_{C\cap S}$$
is injective.
\end{coro}

In the next section this criterion will be applied to conics on prime
Fano threefolds.

\section{Fano surfaces of conics on prime Fano threefolds}

Let $X$ be a prime Fano threefold of genus $g$. It is well known that
$X$ contains lines, but there is only a one dimensional
family of lines on $X$, and in particular $X$ is not covered by lines,
see \S 4.2 in \cite{IP}.
But $X$ is covered by conics, and the Fano scheme $F(X)$ of conics
on $X$ appears as a fundamental object of study, ibid.

In this section we focus on $F(X)$, for $X$ a general prime Fano
threefold, and we apply to conics the main result of the previous
section.

\subsection{Conics on general prime Fano threefolds}

First recall that the genus $g$ of a smooth prime Fano threefold $X$ can
be any integer between $2$ and $12$, except $11$, see \cite{IP}.  

For $g=2$, the
anticanonical divisor is base point free and defines $X$ as
a double cover of $\PP^3$ ramified along a general sextic hypersurface
(a sextic double solid). 
For $g\ge 3$, $-K_X$ is base point free, see Theorem 2.4.5 in \cite{IP}. 
Moreover $-K_X$ is very ample, with only one exception: 
the double quadrics $X_4'$ for which $g = 3$ 
and $-K_X$ defines a double covering from $X$ to the smooth quadric threefold, 
see Proposition 4.1.11 in \cite{IP}. 
Notice that $X_4'$ is not the general prime Fano threefold of genus $g = 3$,
regarded in the context of our main Theorem \ref{main}. 

In all the other cases the smooth $X = X_{2g-2}$ are complete intersections
in homogeneous spaces, 
as follows from the works of Mukai (see \cite{mukai1}, and also
\cite{IP}). 
For the convenience of the reader
we recall below, for each value of $g \ge 3$, the nature of $\Sigma$, its
dimension $d$, and the multidegree $\delta$ of $X$ in $\Sigma$.

$$\begin{array}{llll}
g\hspace*{1cm} & \delta\hspace*{1cm} & \Sigma & d \\
 & & \\
3 & 4&\PP^4& 4 \\
4& 2,3&\PP^5& 5 \\
5&2,2,2 &\PP^6 & 6\\
6&2,1,1 &G(2,5) & 6 \\
7& 1,1,1,1,1,1,1& OG(5,10)& 10\\
8& 1,1,1,1,1&G(2,6)& 8 \\
9& 1,1,1& LG(3,6)& 6\\
10&1,1 & G_2^{ad} & 5
\end{array}$$

\smallskip Here $G(r,d)$ (resp. $OG(r,d)$, $LG(r,d)$) is the Grassmannian
of $r$-dimensional subspaces (resp. isotropic subspaces) of a
$d$-dimensional vector space (resp. endowed with a non degenerate
symmetric of skew-symmetric bilinear form,
where we consider only one connected component
in the case of $OG(r,2r)$). $G_2^{ad}$ denotes the adjoint variety
of the exceptional group $G_2$, the closed $G_2$-orbit in the
projectivized Lie algebra $\PP\fg_2$.

\medskip
For any $2 \le g \le 12, g \not= 11$ the Fano scheme of conics $F(X)$ is known
to be a smooth surface, and the normal bundle to a general conic $q \subset X$
is trivial, see \cite{CV} for $g = 2$, \cite{CMW} for $g = 3$,
\cite{IMn2} for $g = 4$, and Proposition 4.2.5 in \cite{IP} for $g \ge 5$.

\

For $g = 12$, the Fano surface $F(X) \simeq \PP^2$,
see Remark 5.2.15 in \cite{IP}.

\medskip

For $g = 10$, $F(X)$ is an Abelian surface.
Since we could not find any reference for this fact,
we include a proof of the following statement:

\begin{prop}\label{g=10}
For $X$ a general Fano threefold of genus $10$, the Fano surface
$F(X)$ and the intermediate Jacobian $J(X)$ are isomorphic.
\end{prop}

\proof We first prove that $F(X)$ is an Abelian surface.
The idea to prove this is to construct an Abelian fibration on the
Hilbert square of a general K3 surface $S$ of genus $10$.
The existence of this Abelian fibration is known by \cite{HT}
(see also \cite{Saw} for a more general treatment of Abelian
fibrations on Hilbert powers of K3 surfaces).

The general such $S$ is a hyperplane section of the general
Fano threefold $X$ of genus $10$. Following Mukai, $S=\PP^{10}_S\cap
G_2^{ad}$ is a codimension three
linear section of the adjoint variety of $G_2$.
Also the general $X_t$ containing $S$ as a hyperplane section
is a general linear section of $\Sigma$ by a codimension $2$
subspace $\PP^{11}_t$ containing $\PP^{10}$.
Let $\PP^2$ be the family of all $\PP^{10}_t$'s containing $\PP^{10}_S$.
We define  $\cF(S)$ as the family of all conics that lie
on the threefolds $X_t, t \in \PP^2$.
It can be seen that: (a) the 4-fold family $\cF(S)$ is a disjoint union of
the surfaces $F(X_t), t \in \PP^2$ (see Lemma 7 below);
(b) the intersection map
$j:  \cF(S) \rightarrow Hilb_2S, \  q \mapsto q \cap S$
defines an isomorphism between $\cF(S)$ and $Hilb_2S$.
By (a), the family $\cF(S)$ has a natural regular
surjective map $f: \cF(S)  \rightarrow \PP^2$ sending
the conic $q \in \cF(S)$ to the unique $t \in \PP^2$
such that $q \in F(X_t)$.
By (b), the map $f$ defines a regular surjective map
$f_o : Hilb_2S \rightarrow \PP^2$.
By the theorem of Matsushita (see \cite{Msh}),
the general fiber $f_o^{-1}(t)$ of $f_o$
is an Abelian surface in $Hilb_2S$.
Since by construction the general fiber
$f_o^{-1}(t) \simeq f^{-1}(t)$ is isomorphic to $F(X_t)$,
identified with the family of conics on
the general Fano 3-fold $X_t \supset S$
then $F(X_t)$ is an Abelian surface.

\smallskip
The isomorphism of the Abelian surface $F(X)$ with the intermediate
Jacobian of $J(X)$ can be derived on the base of the
birational geometry of $X$ as follows.

Let $\ell \subset X$ be a line. By Theorem 4.3.3(vii) of \cite{IP},
the double projection $\pi_{\ell}$ from $\ell$ defines a birational
isomorphism from $X$ to a 3-fold quadric $Q$.
The double projection $\pi_{\ell}$ contracts the family
of conics intersecting $\ell$ to a smooth curve $\Gamma \subset Q$
of genus $2$ and degree $7$; therefore the intermediate Jacobian
$J(X)$ is isomorphic to the Jacobian $J(\Gamma)$ of $\Gamma$.

By Matsusaka's criterion the Abelian surface $F(X)$ will
be isomorphic to $J(\Gamma)$ if we can find on $F(X)$
a curve $C \simeq \Gamma$  with self-intersection
number $C^2 = 2$, see \cite{Msk}.
We shall show that such a curve is the base $F_{\ell} \subset F(X)$
of the family of conics on $X$ intersecting $\ell$.
First, by the preceding $F_{\ell} \simeq \Gamma$.
Next, the base of the family $G(X)$ of lines on $X$
is a smooth irreducible curve (of genus $10$), see Theorem 4.2.7
in \cite{IP}. Therefore all the curves $F_{\ell}$, $\ell \in G(X)$
are algebraically equivalent to each other on the surface $F(X)$.
So $F_{\ell}^2 = F_{\ell}.F_m$ for any other line $m$
on $X$. If $m \subset X$ is general then the intersection number
$F_{\ell}.F_m$ is the number of conics
$q \subset X$ that intersect both $\ell$ and $m$.
But according to Theorem 4.3.3 (vii) of \cite{IP},
the conics $q \in F_{\ell}$ sweep out a quadratic section
$D_{\ell}$ of $X$. So $D_{\ell}$ intersects the line $m$
at two points, and therefore $F_{\ell}.F_m = D_{\ell}.m = 2$.\qed

\medskip

For $g = 9$, $F(X)$ is a ruled surface with
base a smooth plane quartic $C_X$ such that
$J(X) \simeq J(C_X)$, see \cite{iliev}.

\medskip

For $g = 8$, $F(X)$ is isomorphic to
the Fano surface of lines on the unique smooth
cubic threefold $Y_X$ birational to $X$, \cite{IMn1}.
The Fano surfaces of lines on cubic threefolds
have been studied in deep by Clemens and Griffiths,
\cite{CG}.

\medskip

For $g = 7$, $F(X)$ is isomorphic to the symmetric
square $Sym^2C_X$ of a genus $7$ curve $C_X$,
such that $J(X) = J(C_X)$, see \cite{IMr} or \cite{Kuz}.

\medskip

For $g=6$, the Fano surface $F(X)$ is smooth by \cite{logachev},
Proposition 0.1.

\medskip

For $g=3,4,5$, the Fano threefold $X$ is a complete intersection.
This allows to describe $F(X)$ as the zero-locus of a section of some
globally generated vector bundle; in particular $F(X)$ is smooth for $X$
general, see \cite{CMW} for $g = 3$, \cite{IMn2} for $g=4$, \cite{Wel}
and \cite{PB} for $g=5$.

Finally, the smoothness of $F(X)$ for $g = 2$
is proved by Ceresa and Verra, see Proposition 1.14 in \cite{CV}.

\medskip
We compile some numerical data for the Fano surface $F=F(X)$
of a general prime Fano threefold $X$ of genus $g$.

$$\begin{array}{ccccccc}
 g & c_2 & c_1^2 & p_a & \nu  & q & F \\
 &&&&&& \\
12 & 3 & 9 & 0 & 6 &0& \PP^2 \\
10 & 0 & 0 & -1 & 9  &2& Abelian \\
9 &-8&-16&-3& 12  &3& ruled  \\
8 &27&45&5&16 &5& \\
7 &66&114&14& 24 &7& Sym^2C \\
6 &384&720&91&39&10 \\
5 &1760&3376&427&148&14& \\
4  &11961&23355&2942&318&20& \\
3 &172704&341040&42811&972&30& \\
\end{array}
$$

\medskip
We denoted by $\nu$ the number of conics through a general point of $X$.
For the value of $\nu$ for $ g \ge 8$, see the table (2.8.1) in
\cite{Takeuchi}.
In the cases $g = 7$ and $g = 6$, the numbers $\nu = 24$ and $\nu = 39$
can be computed by the same approach as in \cite{Takeuchi},
see e.g. Theorem 4.5.8 in \cite{IP}. The computation for $g=3$ is done
in \cite{CMW}, and the case $g=4$ is in \cite{IMn2}. The invariants
in genus $5$ are in \cite{Wel}, except for the computation of $\nu$
which can be obtained with the same techniques as for $g=3,4$.

For $g = 2$ (as well for all the other $g$) for the proof of Theorem \ref{main} 
we shall need in fact only the known by \cite{CV} invariant $q = 52$. 
Beware that in the case $g = 6$ the Fano surface $F(X)$ is the blow-up of
a point on a smooth surface $F$; the table gives the invariants of $F$,
see \cite{logachev}.

\subsection{Intersecting conics with a K3 surface}

Let $S\subset X$ be a K3 surface, defined as a general anticanonical
section. Then the image of the natural map $j: F(X)\ra Hilb_2S$
is a possibly singular Lagrangian surface.
In fact $j$ can be an embedding only if $g\ge 6$.
Our precise result is the following:

\begin{prop}\label{immers}
Let $X$ be a general Fano threefold of genus $g$, $S$
a general hyperplane section. Suppose that the curve
$\Gamma(X)$ of lines on $X$ is smooth, and that the
Fano surface $F(X)$ of conics on $X$ is also smooth.
Then:
\begin{enumerate}
\item $j: F(X)\lra Hilb_2S$ is an immersion.
\item $j$ is injective for $g\ge 7$ (but not for $g\le 6$).
\end{enumerate}
\end{prop}

\proof The injectivity of $j$ is equivalent to the non-existence
of bisecant conics meeting along $S$. For $g\ge 7$ this was observed
in \cite{voisin}, Section 2. For $g\le 6$
this is no longer true, see the Remark below.

\smallskip
Nevertheless we can check that $j$ is always immersive by applying
the criterion of Corollary \ref{criterion}.
Let $q\in F(X)$ be a smooth conic. Since $F(X)$ is smooth at $q$
and the normal bundle $N$ of $q$ in $X$ has trivial determinant,
we must have $N=\cO_q\oplus\cO_q$ or
$N=\cO_q(1)\oplus\cO_q(-1)$. In both cases the evaluation map
$$ev_{q\cap S}:H^0(N)\lra N\otimes_{\cO_X}\cO_{q\cap S}$$
is injective, were $q\cap S$ the union of two distinct points
or the tangent direction to $q$ at some point.

Now we consider the case of a singular conic $q=\ell\cup m$, with
vertex $p=\ell\cap m$. Since $q$ is a locally complete intersection,
its normal bundle $N=N_{q/X}=Hom(I_q/I_q^2,\cO_q)$ is a locally
free $\cO_q$-module of rank two. Twisting the exact sequence
$0\ra\cO_q\ra\cO_{\ell}\op\cO_m\ra\cO_p\ra 0$ by $N$, we get
$$0\ra N\ra N(\ell)\op N(m)\ra N_p\ra 0,$$
where $N(\ell)=N\otimes_{\cO_q}\cO_{\ell}$ is a locally-free
$\cO_{\ell}$-module.
Note that since $I_q\supset I_{\ell}I_m$,
$I_mN$ is also a $\cO_{\ell}$-module, again locally free of rank two.
In fact $I_mN=N(\ell)(-p)$, and we have a commutative diagram
$$\begin{array}{ccccccccc}
 & & 0 && 0 &&&& \\
 && \uparrow  && \uparrow &&&& \\
0 & \lra & N_p & \lra & N(\ell)_p\oplus N(m)_p &
\lra &N_p&\lra & 0 \\
 && \uparrow  && \uparrow && || && \\
0&\lra & N&\lra & N(\ell)\op N(m)& \lra & N_p& \lra & 0\\
 && \uparrow  && \uparrow &&&& \\
&& I_mN\oplus I_{\ell}N & = & N(\ell)(-p)\op N(m)(-p)&&&&\\
 && \uparrow  && \uparrow &&&& \\
 & & 0 && 0 &&&&
\end{array}$$

Beware that $I_mN$ is not isomorphic with $N_{\ell}$,
the normal bundle to $\ell$, although these two sheaves are
obviously isomorphic outside $p$. In fact, there is an exact sequence
$$0\lra I_mN \lra N_{\ell}\lra N_{\ell,p}/Im T_pm=\cO_p\ra 0.$$
For $X$ general, the curve $\Gamma(X)$ of lines on $X$ is smooth, hence
$N_{\ell}\simeq\cO_{\ell}\oplus\cO_{\ell}(-1)$. So there are two
possibilities: either $I_mN\simeq\cO_{\ell}(-1)\oplus\cO_{\ell}(-1)$
or $I_mN\simeq\cO_{\ell}\oplus\cO_{\ell}(-2)$, and the latter case
occurs if and only if the map
$$H^0(N_{\ell})\lra N_{\ell,p}/Im T_pm$$
is not injective. Anyway, we get three possibilities:
\begin{enumerate}
\item $N(\ell)$ and $N(m)$ are both trivial, hence $N$ itself is
      trivial;
\item $N_{\ell}$ is trivial but $N(m)=\cO_m(1)\oplus \cO_m(-1)$
(or the opposite): then $H^0(N)\simeq H^0(N(m))$, and each section
is defined on restriction to $\ell$ by its value at $p$;
\item $N(\ell)=\cO_{\ell}(1)\oplus \cO_{\ell}(-1)$ and
$N(m)=\cO_m(1)\oplus \cO_m(-1)$. Note this does not imply that
$N=\cO_q(1)\oplus \cO_q(-1)$: this is  true only if the two
distinguished directions in $N_p$ fit together, which is in
fact excluded by the condition that $h^0(N)=2$.
\end{enumerate}
Now we need to check that
the evaluation map $ev_x\op ev_y : H^0(N)\lra N_x\op N_y$ is
injective, if $x=\ell\cap S$ and $y=m\cap S$ are distinct;
if they coincide, we want $ev^1_x : H^0(N)\lra N\otimes\cO_q/\cM_x^2$
to be injective. For the three possibilities above this is a
straightforward verification. \qed

\begin{coro}
For $X$ a general prime Fano threefold of genus $g\ge 7$,
and $S$ a general hyperplane section, $j(F(X))$ is a smooth
Lagrangian surface in $Hilb_2S$.
\end{coro}

\smallskip\noindent {\bf Remark.}
For $g=6$ recall that $X$ is a complete
intersection of multidegree $\delta=(2,1,1)$ in $G(2,5)$. For any
conic $q$ in $X$, there is a unique hyperplane $H$ in $\CC^5$ such that
$q$ is contained in $G(2,H)$. Since the intersection of $X$ with $G(2,H)$
has degree four, it must be the union of $q$ with an involutive conic
$q'$, which is bisecant to $q$. By an obvious dimension count, there
must exist a finite number of conics $q$ on $X$, whose two common
points with the involutive conic $q'$ lie on $S$ -- and in this case
$j(q)=j(q')$. This shows that $j(F(X))$ is not smooth but has finitely
many double points.

\subsection{Getting back from $j(F(X))$ to $X$}

It is not true in general that $X$ is uniquely determined by
the surface $F(X)$. For example, the Fano surface of the general $X_{22}$
is isomorphic to $\PP^2$, but the Fano 3-folds $X_{22}$ have $6$ moduli,
see e.g. Corollary 1.2 in \cite{Sch}.
Similarly, one can see that the Fano surface of the general $X_{18}$
is the Jacobian of the general genus $2$ curve and has $3$ moduli,
while $X_{18}$ has $10$ moduli.

Here we consider a slightly different question: does
$j(F(X))\subset Hilb_2S$ define $X$ uniquely among those Fano
threefolds that contain $S$?

\begin{theo}\label{back}
For $7\le g\le 10$, a Fano threefold $X$ containing a fixed
K3 surface $S$ as a hyperplane section, is uniquely determined
by the Lagrangian surface $j(F(X))\subset Hilb_2S$.
\end{theo}

We do not know whether this statement is still true for $g\le 6$.
\medskip

\proof We want to recover  $X\supset S$ from the embedded Lagrangian
surface $j(F(X))\subset Hilb_2(S)$.
By \cite{mukai1}, $S$ and  $X$ have linear embeddings
in a rational homogeneous variety $\Sigma=G/P$ of some simple adjoint
Lie group $G$, and this embedding is unique up to the action of $G$.
Therefore we may assume that $S\subset X\subset G/P$.

Now we have to consider separately the different values of $g$.

\medskip \noindent {\bf $g=10$}. Here $G=G_2$ and $G/P$ is the adjoint
$G_2$-variety, the variety parametrizing highest root spaces in
the projectivized Lie algebra $\PP\fg_2$.

\begin{lemm}
Consider two general points $x,y$ in the adjoint variety $G_2/P$.
Then there is a unique conic in $G_2/P$ passing through $x$ and $y$.
\end{lemm}

\noindent {\it Proof of the lemma}.
First, there exists at least one conic in $G_2/P$ through $x$ and $y$.
Indeed, the Lie subalgebra of $\fg_2$ generated by $x$ and $y$ is a
copy of $\fsl_2$, and generates in $G_2$ a copy of $SL_2$. The
$SL_2$-orbits of $x$ and $y$ have the same closure, which is a conic
passing through $x$ and $y$ in $G_2/P$. (To check this, choose
a maximal torus in $\fg$ with respect to which $x$ and $y$ are
highest and lowest root vectors.)

On the other hand,
the group $G_2$, being the automorphism group of the (complexified)
algebra of octonions, has a minimal representation of dimension $7$
on the imaginary octonions. Since it preserves the norm, we deduce in
particular that the adjoint variety $G_2/P$ is a sub-variety of
$G_Q(2,7)$, the Grassmannian of lines in a smooth quadric
$Q \subset P^6$. But through two general points of a Grassmannian
of lines in a smooth quadric there exists a unique conic, which can be
constructed as follows. The two points $x$ and $y$ represent two skew lines
$\ell_x, \ell_y \subset Q$. These two lines generate
$P^3_{x,y} = Span(\ell_x \cup\ell_y)$ intersecting $Q$ along
a quadric surface $S_{x,y}$. One of the
two rulings of this surface contains $\ell_x,\ell_y$;
it defines a conic $q$ in $G_Q(2,7)$ passing through $x$ and $y$.
\qed

\medskip
From $j(F(X)) \subset Hilb_2S$ we can therefore recover
the Fano surface $F(X)$, and then the threefold $X$ itself
since it is swept out by the conics it contains.

\medskip \noindent {\bf $g=9 $}.
Here $G=Sp_6$ and $G/P=LG(3,6)$, the Lagrangian Grassmannian.

\smallskip
Recall from \cite{iliev,IR} that a conic $q \subset LG(3,6) =
LG(\PP^2,\PP^5)$ has a vertex-point
$v(q) \in \PP^5$, defined as the unique intersection point
of the Lagrangian planes $\PP^2_t$, $t \in q$.

A general conic $q \subset X$ is smooth and intersects
$S$ at two different points $x,y$, from which the vertex can
be recovered  as $v(q) = \PP^2_x \cap \PP^2_y$.
In particular, $j(F(X))$ determines the vertex surface $v(F(X))$.
By \cite{iliev}, for a general Fano threefold $X$ of genus $8$
the vertex image $v(F(X)) \subset \PP^5$ is a ruled surface
isomorphic to $F(X)$. By \cite{IR}, Proposition 2.8.5, $X$ can
be recovered as the set of projective planes in $\PP^5$ that cut
the vertex surface along a finite  scheme of length 12.

\medskip \noindent {\bf $g=8 $}.
Here $G=PSL_6$ and $G/P=G(2,6)$ is the Grassmannian of (projective)
lines.

\smallskip If $S\subset X\subset G(2,6)\subset\PP(\wedge^2\CC^6)$,
let $L_X\subset L_S$ be the orthogonal in $\PP(\wedge^2\CC^6)^*$ of
the linear spans of $S$ and $X$. They are projective spaces of dimension
four and five respectively. We want to show that, the embedding
$j(F(X))\subset Hilb_2S$ being given, the hyperplane $L_X$ of  $L_S$
is uniquely determined.

A conic $q$ in $G(2,6)$ whose linear span is not contained in $G(2,6)$,
defines a unique 4-space $V\subset\CC^6$ such that $q\subset G(2,V)$.
Moreover, $q$ can be obtained as the intersection of the four
dimensional quadric $G(2,V)\subset\PP\wedge^2V$ with some projective
plane. In particular, $G(2,V)$ contains some conic of $X$ if and
only if $\PP\wedge^2V$ meets $L_X^{\perp}$ along a projective space of
dimension two at least. We thus consider on $G=G(4,6)$ the tautological
rank four vector bundle $T$. By restriction of skew-symmetric forms,
we have morphisms  of vector bundles
$$\begin{array}{ccc}
L_S\otimes \cO_G & \stackrel{\phi_S}{\lra} & \wedge^2T^* \\
\uparrow &  & || \\
L_X\otimes \cO_G & \stackrel{\phi_X}{\lra} & \wedge^2T^*
\end{array}$$
Since $\phi_X$ is a general morphism from $\cO_G^{\oplus 5}$
to the rank six vector bundle $\wedge^2T^*$, its degeneracy locus
$\Sigma(X)=D_3(\phi_X)$ has codimension $(5-3)(6-3)=6$ and is singular
exactly along $D_2(\phi_X)$, which has codimension $(5-2)(6-2)=12>8$.
So $\Sigma(X)$ is a smooth surface, and since it is in bijective
correspondence with $F(X)$ these two surfaces must be isomorphic.

On the other hand, the degeneracy locus
$\Sigma(S)=D_4(\phi_S)$ has codimension $(6-4)(6-4)=4$ and is singular
exactly along $D_3(\phi_S)$, which has codimension $(6-3)(6-3)=9>8$.
So $\Sigma(S)$ is a smooth four dimensional subvariety of $G(4,6)$.
On $\Sigma(S)$ the morphism $\phi_S$ has constant rank four, so in the
exact sequence
$$0\lra Ker\phi_S\lra L_S\otimes \cO_G  \stackrel{\phi_S}{\lra}
\wedge^2T^*\lra Coker\phi_S\lra 0,$$
$Ker\phi_S$ and $Coker\phi_S$ are rank two vector bundles.
Note that if we take determinants in this sequence, we get
$$\det(Coker\phi_S)\det(Ker\phi_S)^{-1}=\det(\wedge^2T^*)=\cO_G(3).$$
Now the normal bundle of $\Sigma(S)$ in $G$ is naturally identified
with $Hom(Ker\phi_S,Coker\phi_S)$, whose determinant is
$\det(Coker\phi_S)^2\det(Ker\phi_S)^{-2}=\cO_G(6)=K_G^{-1}$. Therefore
$\Sigma(S)$ has trivial canonical bundle. Now the natural map
$\Sigma(S)\ra Hilb_2S$ is finite, hence \'etale, hence an isomorphism
since $Hilb_2S$ is simply connected. We conclude that the embedding
$j(F(X))\subset Hilb_2S$ is the same as the embedding $\Sigma(X)\subset
\Sigma(S)$.

But from $L_S$, its degeneracy locus $\Sigma(S)$ and the embedding
$\Sigma(X)\subset \Sigma(S)$, the hyperplane $L_X$ of $L_S$ can be
recovered as follows. At any point $x\in\Sigma(X)$, the morphisms
$\phi_S$ ans $\phi_X$ have the same kernel, which in particular
is contained in $L_X$. So we must have
$$L_X=Span(Ker\phi_{S,x}, \; x\in \Sigma(X))\subset L_S.$$
And then of course $X$ can be recovered as $X=G(2,6)\cap L_X^{\perp}.$

\medskip \noindent {\bf $g=7 $}.
Here $G=Spin_{10}$ and $G/P=\Sigma\subset\PP\Delta_+$ is the spinor
variety, naturally embedded in the projectivization of a half-spin
representation.

\smallskip If $S\subset X\subset\Sigma\subset\PP\Delta_+$,
let $L_X\subset L_S$ be the orthogonal in $\PP\Delta_-$ of
the linear spans of $S$ and $X$. They are projective spaces of dimension
seven and eight respectively. Again we want to show that, the embedding
$j(F(X))\subset Hilb_2S$ being given, the hyperplane $L_X$ of  $L_S$
is uniquely determined.

A conic $q$ in $\Sigma$ whose linear span is not contained in $\Sigma$,
has a vertex $v=v(q)\in\PP^9$. The set of isotropic planes in $\Sigma$
passing through $v$ is a copy of the spinor variety of $Spin_8$. By
triality, this spinor variety is a six dimensional quadric $Q_v\subset
\Sigma$. Moreover, $q$ can be obtained as the intersection of $Q_v$
with some projective plane. In particular, $Q_v$ contains some conic
of $X$ if and only if its linear space $P_v$ meets $L_X^{\perp}$ along
a projective space of dimension two at least.

This linear space $P_v$ defines on the eight-dimensional quadric $\QQ^8$
in $\PP^9$ a spinor bundle $S$ of rank eight \cite{ott}. We get on $\QQ^8$
the morphisms  of vector bundles
$$\begin{array}{ccc}
L_S\otimes \cO_G & \stackrel{\phi_S}{\lra} & S^* \\
\uparrow &  & || \\
L_X\otimes \cO_G & \stackrel{\phi_X}{\lra} & S^*
\end{array}$$
Exactly as in genus eight, the degeneracy locus $\Sigma(X)=D_3(\phi_X)$
 is a smooth surface isomorphic with $F(X)$. Also the degeneracy locus
$\Sigma(S)=D_2(\phi_S)$ is a smooth four dimensional variety with
trivial canonical bundle (because $\det S=\cO_Q(4)$), in fact isomorphic
with $Hilb_2S$. And from the embedding $\Sigma(X)\subset \Sigma(S)$,
which is the same as $j(F(X))\subset Hilb_2S$, we can recover $L_X$
in $L_S$ as before.

\smallskip
This concludes the proof. \qed

\section{Fano-Integrable Systems}

Thanks to the main results of the previous section, we are now in
position to apply the Donagi-Markman construction \cite{DM}.
The upshot will be new Lagrangian fibrations over families of prime
Fano threefolds with a fixed ``boundary'', whose fibers are intermediate
Jacobians.

\subsection{The Abel-Jacobi isomorphism for conics on prime Fano's}

The intermediate Jacobian $J(X)$ a prime Fano threefold of genus
$g$ is an Abelian variety of dimension $h^{1,2}(X)=b_3(X)/2$ (see
\cite{CG} for the definition and general properties of intermediate
Jacobians).

This dimension can be computed case by case.
Note that for $7\le g\le 10$, we can use the fact that $X$ is a linear
section of a $d$-dimensional homogeneous space
$\Sigma\subset \PP^{d+g-2}$. After twisting the normal exact sequence
by $K_X=\cO_X(-1),$ we get
$$0\lra \Omega_X^2\lra T\Sigma\otimes\cO_X(-1)\lra
\cO_X^{\oplus d-3}\lra 0.$$
With the help of the Bott-Borel-Weil theorem one can check
that $H^{i}(\Sigma,T\Sigma\otimes\cO_{\Sigma}(-i-1))=
H^{i+1}(\Sigma,T\Sigma\otimes\cO_{\Sigma}(-i-1))=0$ for
$0\le i\le d-3$. Using the Koszul resolution of $\cO_X$
this implies that $H^0(X,T\Sigma\otimes\cO_X(-1))=
H^1(X,T\Sigma\otimes\cO_X(-1))=0$. Hence the isomorphism
$$H^1(\Omega_X^2)=\CC^{d-3}=L_X,$$
if we denote by $L_X$ the space of linear forms defining $X$ in
$\Sigma$.

Note also that if $X$ has no automorphism, which is the general case,
we get an exact sequence
$$0\lra H^0(X,T\Sigma_{|X})\simeq\fg\lra L_X\otimes L_X^{\perp}
\lra H^1(X,TX)\lra 0,$$
where $\fg$ denotes the Lie algebra of the Lie group $G$;
hence $h^1(TX)=(g+2)(d-3)-\dim G$. The analogue sequence for $S$ is
$$0\lra H^0(S,T\Sigma_{|S})\simeq\fg\lra L_S\otimes L_S^{\perp}
\lra H^1(S,TS)\lra \CC\lra 0.$$
This gives $h^1(TS)-1=(g+1)(d-2)-\dim G=19$, and thus
$$h^1(TX)=15+d-g=h^1(\Omega_X^2)+18-g.$$

\smallskip
A basic property of the intermediate Jacobian is that for any
family $\cC\ra F$ of curves in $X$, parameterized by a smooth base
$F$, there is an induced {\it Abel-Jacobi mapping} $F\ra J(X)$.

\begin{theo}\label{AJ}
Let $X$ be a general prime Fano threefold of genus $g$, with a smooth
 surface of conics $F(X)$.
Then the Abel-Jacobi mapping $F(X)\ra J(X)$ induces an
isomorphism
$$Alb F(X)\simeq J(X).$$
\end{theo}

\proof 
For most values of $g$, this Abel-Jacobi isomorphism is already 
directly or indirectly known, so below we refer to the original 
references, interpret some other in the required context, and 
complete the argument for the genera $g = 10$ and $6$ missing 
in the known to us literature.

The case $g=12$ is trivial since then $h^{1,2}(X)=0$ 
and $F(X)$ is isomorphic to $\PP^2$, see \S 3.1. 

For $g = 10$,  by Proposition \ref{g=10} the Fano surface 
$F(X)$ is isomorphic to $J(X)$ so in this case there 
is already nothing more to prove. 

For $g=9$, $F(X)$ is a ruled surface over the orthogonal plane
quartic curve $C_X$, therefore $Alb F(X)=J(C_X)$, see \S 3.1. 
The fact that $J(C_X)\simeq J(X)$ is proved in \cite{IR}, 
Corollary 2.8.10.

For $g=8$, $F(X)$ is isomorphic to the Fano surface $F(Y)$
of lines in the orthogonal cubic threefold $Y$, see \S 3.1. 
The intermediate
Jacobians of $X$ and $Y$ are naturally isomorphic. The Fano
surface $F(X)$ embeds in $J(X)$, $F(Y)$ embeds in $J(Y)$, and
their images coincide up to a multiplication by $-1$ (\cite{IMn1},
Proposition 4.6). So the Abel-Jacobi isomorphism for $F(X)$
follows from the corresponding statement for $F(Y)$, which is
due to Clemens and Griffiths \cite{CG}.

For $g=7$, $F(X)$ is isomorphic to $Sym^2C_X$, where $C_X$ is
the $Spin_{10}$-orthogonal curve to $X$, see \S 3.1. From where
$Alb F(X)\simeq J(C_X)$. As known by \S 8 of \cite{mukai2}
the Jacobian $J(C_X)$ is isomorphic to $J(X)$,
see also \cite{IMr}.

For $g = 3$, the Abel-Jacobi isomorphism for $F(X)$ 
has been proved by Letizia, by applying a degeneration technique 
originating from Clemens, see \cite{Let}. 
After that, the same method had been used by Ceresa and Verra
(see \cite{CV}), Picco Botta (see \cite{PB}) 
and in our recent paper \cite{IMn2} to prove the same result 
correspondingly for $g = 2,5$ and $4$. 

The case $g = 6$ has been treated by Logachev
(before the appearance of the paper \cite{Let})   
with the conclusion that $Alb F(X)\ra J(X)$ in an isogeny, 
by following an earlier argument of Welters 
proving the same isogeny for the case $g = 5$, 
see Theorem 6.12 in \cite{logachev} and \cite{Wel}. 
In fact, the Clemens-Letizia method applies 
also in the case $g = 6$ to conclude that the 
Abel-Jacobi map is an isomorphism also for $g = 6$. 
All the necessary conditions for this 
(a precise description of the Fano surface for a general nodal 
Fano threefold) are provided in Section 5 of \cite{logachev}.
\qed

\subsection{The intermediate Jacobian integrable systems for prime Fano's}

Consider a projective symplectic  variety $Y$ of dimension $2n$,
a smooth variety $L$ of dimension $n$ and a Lagrangian immersion
$f: L\ra Y$. Then:
\begin{enumerate}
\item If $\cL\ra\cH$ is a local semi-universal deformation of  $f$, with fixed
      target $Y$, then $\cH=Def(f/Y)$ is smooth (\cite{Ran}, Corollary
      3.4).
\item The tangent space to $\cH$ at $L$ is
$$T_{f}\cH=H^0(L,N_f)\simeq H^0(L,\Omega_L^1),$$
where the bundle $N_f$ is defined by the exact sequence
$0\ra TL\ra f^*TY\ra N_f\ra 0$; indeed
the fact that $f$ is Lagrangian yields an isomorphism
$N_f\simeq \Omega_L^1$.
\item Any small deformation of $f$ remains a Lagrangian immersion in $Y$,
see \cite{voisin}.
\item The relative Picard bundle $Pic^0\cL\lra\cH$ is a Lagrangian
      fibration  (or an algebraic
completely integrable Hamiltonian system), see
\cite{DM}. In fact Donagi and Markman only treat the case where
$j$ is an embedding, that is, they consider the deformations of
a smooth Lagrangian subvariety in a symplectic manifold. Nevertheless,
their proof extends verbatim to the setting of Lagrangian immersions.
\end{enumerate}

Let $X$ be a general prime Fano threefold of genus $g$ and $S$
a K3 surface in $X$. We can apply the results above to $f=j: L=F(X)
\ra Y=Hilb_2S$, since we have proved in Proposition \ref{immers}
that $j$ is indeed a Lagrangian immersion. We get that the family
$Pic^0\cL$ over $Def(j/Hilb_2S)$ is a Lagrangian fibration.

\smallskip
Now we want to deform $X$ with $S$ fixed. The tangent space to the
corresponding deformation functor is $H^1(X,TX\otimes\cI_S)$
(see \cite{beau3}). Since $\cI_S=\cO_X(-S)=K_X$, this is also
$H^1(\Omega^2_X)$, the universal covering of $J(X)$.

We begin with a simple observation.

\begin{prop}
There exists a semi-universal algebraic family $\cH_S$
parameterizing deformations of $X$ containing $S$.
\end{prop}

\proof Suppose for example that $7\le g\le 10$ (in which case $j$ is
an embedding). Recall that by Mukai's theorem, $S$ can be described as
a linear section $S=\Sigma\cap L_S^{\perp}$ of a $d$-dimensional
homogeneous space $\Sigma$ by a projective space of codimension $d-2$.
The set of prime Fano threefolds $X$ containing $Š$ as an anticanonical
section, is then identified with an open subset of the orthogonal
projective space $L_S\simeq\PP^{d-3}$.

In the complete intersection case, consider for example the general
quartic threefold $X\subset\PP^4$ admitting $S$ as a hyperplane section.
The set of such $X$'s is given by the set of quartic polynomials
of the form $k=k_S+x_0c$, where $x_0$ is an equation of the hyperplane
$H$ spanned by $S$ in $\PP^4$, $k_S$ is an equation of $S$ in $H$, and
$c$ is any cubic form. So the smooth quartic hypersurfaces of $\PP^4$
containing $S$ are parameterized by an open subset of a $35$-dimensional
affine space. Moreover the $5$-dimensional subgroup $\Gamma_H$ of
$GL_5$, acting trivially on $H$, acts on this affine space, and any
linear subspace transverse to the orbit of $X$ defines a semi-universal,
$30$-dimensional, deformation of $X$ with $S$ fixed. The other cases
are similar. \qed

\medskip
We have a natural map from the germ $(\cH_S,X)$ to $Def(j/Hilb_2S)$,
 whose differential
$$H^1(X,TX(-S))=H^1(X,\Omega_X^2)\lra H^0(F(X),\Omega_{F(X)}^1)$$
is an isomorphism when the Abel-Jacobi isomorphism does hold.
In particular we can pull-back our Lagrangian fibration
to $\cH_S$, first locally, and then everything globalizes since the
fibration is defined canonically -- more precisely, we can pull it back
to the open subset $\cH_S^0$ on which $X$ is smooth, $F(X)$ is also
smooth and the Abel-Jacobi isomorphism theorem does hold.
But then over $\cH_S^0$, $Pic^0F(X)\simeq Alb F(X)^*
\simeq J(X)^*\simeq J(X),$
the intermediate Jacobian being self-dual. Since this isomorphism
certainly holds in families, we finally deduce our main result:

\begin{theo}\label{main} 
Let $S$ be a general K3 surface of genus $g$,
with $2\le g\le 10$. Let $\cX\ra\cH_S$ be an algebraic semi-universal family
of smooth Fano threefolds containing $S$.

Then the family of intermediate
Jacobians $J(\cX)\ra \cH_S$ is a Lagrangian fibration over a dense
open subset $\cH_S^0$ of $\cH_S$.
\end{theo}

For $g\ge 7$ the situation is even cleaner: one can take for $\cH_S$
an open subset of the projective space $L_S$, and
by Theorem \ref{back} there is a birational isomorphism between
$L_S$ and the component $\cH$ of the Hilbert
scheme of $Hilb_2S$ parameterizing the Lagrangian
surfaces $jF(X)$, for $X$
a general Fano threefold containing $S$.


\vspace{1cm}

{\small
\noindent
{\bf Atanas Iliev}\\
Institute of Mathematics\\
Bulgarian Academy of Sciences\\
Acad. G. Bonchev Str., bl.8\\
1113 Sofia, Bulgaria\\
{\bf e-mail:} ailiev@math.bas.bg

\bigskip

\noindent
{\bf Laurent Manivel}\\
Institut Fourier\\
UMR 5582 (UJF-CNRS), BP 74\\
38402 St Martin d'H\`eres Cedex, France\\
{\bf e-mail:}  Laurent.Manivel@ujf-grenoble.fr
}

\end{document}